\documentclass[11pt]{article}%
\usepackage{amsmath}
\usepackage{amsfonts}
\usepackage{amssymb}
\usepackage{graphicx}
\usepackage{theorem}
\usepackage{color}%
\providecommand{\U}[1]{\protect\rule{.1in}{.1in}}
\newtheorem{thm}{Theorem}[section]
\newtheorem{lm}[thm]{Lemma}
\newtheorem{pr}[thm]{Proposition}
\newtheorem{df}[thm]{Definition}
\newtheorem{rmk}[thm]{Remark}
\newtheorem{cor}[thm]{Corollary}
{\theorembodyfont{\upshape}
\newtheorem{examp}[thm]{Example}
}
\numberwithin{equation}{section} 
\begin{document}

\title{Title: }

\begin{center}
\bigskip

\vspace*{1.3cm}

\textbf{CONIC CANCELLATION LAWS AND SOME\ APPLICATIONS}

\bigskip

by

\bigskip

Marius\ DUREA\footnote{{\small Faculty of Mathematics, \textquotedblleft
Alexandru Ioan Cuza\textquotedblright\ University, 700506--Ia\c{s}i, Romania
and \textquotedblleft Octav Mayer\textquotedblright\ Institute of Mathematics,
Ia\c{s}i Branch of Romanian Academy, 700505--Ia\c{s}i, Romania{; e-mail:
\texttt{durea@uaic.ro}}}} and Elena-Andreea FLOREA\footnote{{\small Faculty of
Mathematics, \textquotedblleft Alexandru Ioan Cuza\textquotedblright%
\ University, 700506--Ia\c{s}i, Romania and \textquotedblleft Octav
Mayer\textquotedblright\ Institute of Mathematics, Ia\c{s}i Branch of Romanian
Academy, 700505--Ia\c{s}i, Romania{; e-mail:
\texttt{andreea\_acsinte@yahoo.com}}}}
\end{center}

\bigskip

\noindent{\small {\textbf{Abstract:}} We discuss, on finite and infinite
dimensional normed vector spaces, some versions of R\aa dstr\"{o}m
cancellation law (or lemma) that are suited for applications to set
optimization problems. In this sense, we call our results "conic" variants of
the celebrated result of R\aa dstr\"{o}m, since they involve the presence of
an ordering cone on the underlying space. Several adaptations to this context
of some topological properties of sets are studied and some applications to
subdifferential calculus associated to set-valued maps and to necessary
optimality conditions for constrained set optimization problems are given.
Finally, a stability problem is considered.}

\bigskip


\bigskip


\begin{center}

\end{center}

\section{Introduction and notation}

In this paper we develop some ideas firstly emphasized in \cite{DF2023-Subd}
where subdifferential calculus rules for set-valued maps were presented.
Actually, in turn, the generalized subgradients under study in that paper
(primarily defined in \cite{DS23-JOGO}) are designed to deal with set
optimization problems and therefore are based on epigraphical associated
set-valued maps, where the epigraphs are defined by means of the ordering cone
that shapes the optimization problems. It was apparent from that study that
for some such calculus rules, a version of R\aa dstr\"{o}m cancellation law
involving the presence of the ordering cone is useful and this gives us the
impetus to further explore the subject. Consequently, in this paper we study
several conic variants of cancellation laws on infinite and finite dimensional
settings and we present some applications in the topic which started this
investigation, namely, calculus for subdifferentials of set-valued maps.

We base our investigation in the main section of this work (that is, Section
2) on several tools among which we mention linear and nonlinear separation
results, weak and strong compactness of a set with respect to a cone and a
metric regularity condition for sets. In this sense, we use classical
separation results for convex sets and also the Gerstewitz (Tammer)
scalarizing functional which provides nonconvex separation under certain
conditions. Moreover, besides compactness of a set with respect to a cone in
the standard meaning (see \cite{Luc}, \cite{DF2022}) we propose a weaker
notion which considers the weak topology of the underlying space. Furthermore,
on finite dimensional spaces we make use of some regularity properties of sets
investigated in \cite{BB} which evaluate the distance from a point to the
intersection of two sets by the sum of the distance from that point to each
individual set.

The third section uses some of the cancellation rules presented before in
order to derive results concerning the invariance of the excess (and,
implicitly, of the Hausdorff distance) to the addition of a set in both terms
and to present calculus rules for the generalized subgradients involving
set-valued maps. Several comments that show the possibility to use the
embedding approach discussed in \cite{KN}, \cite{SK} (see also \cite{BT}) for
dealing with set optimization in a broader context, using a class of unbounded
sets, are presented. Finally, we introduce a concept of sharp minimality for
constrained set optimization problems and we present a necessary optimality
condition for whose proof we employ the R\aa dstr\"{o}m cancellation law.
Moreover, we show how a stability principle devised in \cite{Sha} for scalar
optimization problems can be extended in the current framework by using a
conic cancellation law.

\bigskip

The notation is standard. Let $X$ be a normed space over the real field
$\mathbb{R}$. The topological dual of $X$ is $X^{\ast},$ the norm will be
denoted $\left\Vert \cdot\right\Vert .$ If $Z$ is also a normed vector space
we denote by $B\left(  Z,X\right)  $ the normed vector space of linear bounded
operators from $Z$ to $X.$ We put $B\left(  x,\varepsilon\right)  $ for the
open ball centered at $x\in X$ with the radius $\varepsilon>0,$ while we
denote by $B_{X},$ $D_{X}\ $and $S_{X}$ the unit open ball, the unit closed
ball and the unit sphere of $X,$ respectively. If $A\subset X$ is a nonempty
set, $\operatorname{cl}A,$ $\operatorname*{int}A,$ $\operatorname*{conv}A,$
$d\left(  \cdot,A\right)  $ are the topological closure, the topological
interior, the convex hull, and the associated distance function, respectively.
If $A,B$ are nonempty subset of $X$, the excess from $A$ to $B$ is
\[
e(A,B)=\sup_{x\in A}d(x,B).
\]
It is easy to see that
\begin{equation}
e(A,B)=\inf\left\{  \alpha>0\mid A\subset B+\alpha B_{X}\right\}
=\inf\left\{  \alpha>0\mid A\subset B+\alpha D_{X}\right\}  , \label{rel1}%
\end{equation}
and
\begin{equation}
e(A,B)=e(\operatorname{cl}A,B)=e(A,\operatorname{cl}B)=e(\operatorname{cl}%
A,\operatorname{cl}B), \label{rel2}%
\end{equation}
where the usual convention $\inf\emptyset=+\infty$ is in use.

\section{Conic cancellation laws}

The R\aa dstr\"{o}m cancellation law is a well-known result (see \cite[Lemma
1]{Rad} and \cite[Lemma 2.1]{Sch}).

\begin{pr}
Let $A,B,C$ be given nonempty sets in $X$ and suppose that $B$ is closed and
convex, $C$ is bounded, and $A+C\subset C+B.$ Then $A\subset B.$
\end{pr}

In order to deal with conic versions of this result we have to recall some
fundamental notions. We consider $K\subset X$ to be a closed convex and
pointed cone (in the sequel, this is the standing assumption for $K$). We
denote by $K^{+}$ its positive dual cone and we say that $K$ is solid if
$\operatorname*{int}K\neq\emptyset.$ In \cite{DF2022} the notion of
$K-$sequential compactness was introduced and studied, namely, a nonempty set
$A$ is called $K-$sequentially compact if it satisfies the following property:
for any sequence $\left(  a_{n}\right)  \subset A$ there is a sequence
$\left(  c_{n}\right)  \subset K$ such that the sequence $\left(  a_{n}%
-c_{n}\right)  $ has a convergent subsequence (in norm topology) towards an
element in $A$. Recall as well (see \cite{Luc}) that the nonempty set $A$ is
called $K-$bounded if there is a bounded set $M\subset X$ such that $A\subset
M+K$, $K-$closed if $A+K$ is closed, and $K-$convex if $A+K$ is convex. Notice
that a $K-$sequentially compact set is $K-$bounded and $K-$closed (see
\cite{DF2022}).

In \cite{DF2023-Subd}, a "conic" cancellation law was proved. We reproduce it
in a slightly weaker form which is suited for the subsequent discussion. For
completeness, and also in order to mark the difference in conception with the
new cancellation laws we are going to investigate in this work, we give a
sketch of the proof.

\begin{pr}
\label{pr_Rad_conic}Suppose that $A,B,C\subset X\ $are nonempty sets such that
$C$ is $K-$bounded and
\[
A+C\subset C+B+K.
\]
Then
\[
A\subset\operatorname{cl}\operatorname*{conv}\left(  B+K\right)  .
\]

\end{pr}

\noindent\textbf{Proof. }Observe that it is enough to take $A$ as a singleton,
$A=\left\{  a\right\}  $ with $a\in X.$ Since $\left\{  a\right\}  +C\subset
C+B+K,$ one has $C\subset C+\left(  B-a\right)  +K$ and
\[
\operatorname*{cl}\operatorname*{conv}\left(  \left(  B+K-a\right)  \right)
=\operatorname*{cl}\left(  \operatorname*{conv}\left(  B+K\right)  -a\right)
=\operatorname{cl}\operatorname*{conv}\left(  B+K\right)  -a,
\]
we can even consider $A=\left\{  0\right\}  .$ Denote by $M$ a bounded set
that satisfies $C\subset M+K.$ Then for all natural $n,$%
\begin{align*}
C  &  \subset B+K+C\subset B+K+B+K+C\subset...\\
&  \subset\underset{n\text{ times}}{\underbrace{B+...+B}}+K+C\subset
\underset{n\text{ times}}{\underbrace{B+...+B}}+K+M.
\end{align*}
Take $c\in C.$ Then for all $n$ there is $m_{n}\in M$ such that%
\[
\frac{1}{n}\left(  c-m_{n}\right)  \in\frac{1}{n}\left(  \underset{n\text{
times}}{\underbrace{B+K+...+B+K}}\right)  \subset\operatorname*{conv}\left(
B+K\right)  .
\]
Since $M$ is bounded, we get that $0\in\operatorname{cl}\operatorname*{conv}%
\left(  B+K\right)  $ and this is the conclusion.\hfill$\square$

\bigskip

In this section we explore more conic variants of this result in different
settings. Firstly, we are interested in other variants of R\aa dstr\"{o}m
cancellation law on infinite dimensional spaces. The prototype is going to be
the following (non-conic) cancellation law. One can see that the proof of this
result is based on a separation result (for convex sets) and this idea is
constantly used for the rest of the section.

\begin{pr}
Suppose that $A,B,C\subset X\ $are nonempty sets, such that $C$ is weakly
compact, $B$ is open and
\[
A+C\subset C+B.
\]
Then
\[
A\subset\operatorname*{conv}B.
\]

\end{pr}

\noindent\textbf{Proof. }As above, it is enough to take $A$ as a singleton and
by a translation we can even consider $A=\left\{  0\right\}  .$ Suppose that
$0\notin\operatorname*{conv}B.$ Then, by a separation result, since $B$ is
open and the convex hull of an open set is open, there is $x^{\ast}\in
X^{\ast}$ such that%
\[
0<x^{\ast}\left(  b\right)  ,\text{ }\forall b\in\operatorname*{conv}B.
\]
By the weak compactness of $C,$ the functional $x^{\ast}$ achieves its minimum
of $C,$ so there is $\overline{c}\in C$ such that
\[
x^{\ast}\left(  \overline{c}\right)  \leq x^{\ast}\left(  c\right)  ,\text{
}\forall c\in C.
\]
Then there are $c\in C,$ $b\in B$ such that $\overline{c}=c+b,$ so%
\[
x^{\ast}\left(  \overline{c}\right)  =x^{\ast}\left(  c\right)  +x^{\ast
}\left(  b\right)  \geq x^{\ast}\left(  \overline{c}\right)  +x^{\ast}\left(
b\right)  >x^{\ast}\left(  \overline{c}\right)  .
\]
This is a contradiction, so the conclusion holds.\hfill$\square$

\bigskip

Now, in order to get conic variants of this result, we have to prepare some
appropriate tools. We recall (see \cite{Luc}) that a nonempty subset $A\subset
X$ is called $K-$compact (or compact with respect to the cone $K$) if from any
cover of $A$ with the sets of the form $U+K$, where $U$ is open (in norm
topology), one can extract a finite subcover of it. It is shown in
\cite{DF2022} that such a set is also $K-$sequentially compact, and, moreover,
the converse holds provided $K\ $is separable.

We work in the sequel with a weaker property, which we call weakly
$K-$compactness and which is defined by simply taking the sets $U$ in the
definition of $K-$compactness as being weakly open.

\begin{examp}
Take $X=\ell^{2}$ with its usual norm and $K=\ell_{+}^{2}.$ Then the closed
unit ball is weakly $K-$compact and it is not $K-$compact. Indeed, the first
assertion is obvious since $D_{\ell^{2}}$ is weakly compact. For the second
assertion, taking into account the above comments, it is enough to prove that
$D_{\ell^{2}}$ is not $K-$sequentially compact. For this, take the sequence
$\left(  -e_{n}\right)  _{n\geq1}\subset D_{\ell^{2}}$, where $\left(
e_{n}\right)  _{n\geq1}$ denotes the sequence of the unit vectors. Suppose
that there is $\left(  c_{n}\right)  _{n\geq1}\subset\ell_{+}^{2}$ such that
$\left(  -e_{n}-c_{n}\right)  _{n}$ has a subsequence (denoted the same)
convergent in norm topology towards an element in $D_{\ell^{2}}.$ This means
that for all $\varepsilon>0$ there is $n_{\varepsilon}\in\mathbb{N}$ such that
for all $n\geq n_{\varepsilon},$%
\[
\sum_{p=1}^{\infty}\left(  -e_{n}^{p}-c_{n}^{p}\right)  ^{2}<1+\varepsilon,
\]
where the superscript is numbering the terms of the underlying element
(sequence) of $\ell^{2}.$ Taking into account the particularities of unit
vectors and that $\left(  c_{n}\right)  \subset\ell_{+}^{2}$, we get that
\[
\sum_{p=1}^{\infty}\left(  c_{n}^{p}\right)  ^{2}<\varepsilon,\text{ }\forall
n\geq n_{\varepsilon}.
\]
This shows that $\left(  c_{n}\right)  $ has to be convergent to $0,$ whence,
in particular, a Cauchy sequence. Consequently, for all distinct and large
enough $n,m\in\mathbb{N}$ we have%
\begin{align*}
\left\Vert \left(  -e_{n}-c_{n}\right)  -\left(  -e_{m}-c_{m}\right)
\right\Vert _{2}  &  =\left\Vert \left(  e_{m}-e_{n}\right)  -\left(
c_{n}-c_{m}\right)  \right\Vert _{2}\geq\left\Vert e_{m}-e_{n}\right\Vert
_{2}-\left\Vert c_{n}-c_{m}\right\Vert _{2}\\
&  =\sqrt{2}-\left\Vert c_{n}-c_{m}\right\Vert _{2}>1.
\end{align*}
We conclude, that, actually, $\left(  -e_{n}-c_{n}\right)  _{n}$ is not a
Cauchy sequence, and this is a contradiction.
\end{examp}

\begin{lm}
\label{lm_min}If $x^{\ast}\in K^{+}$ and $A$ is weakly $K-$compact, then
$x^{\ast}$ achieves its minimum on $A.$
\end{lm}

\noindent\textbf{Proof. }Denote $\alpha=\inf\left\{  x^{\ast}\left(  a\right)
\mid a\in A\right\}  .$ Suppose, by way of contradiction, that $\alpha
<x^{\ast}\left(  a\right)  $ for all $a\in A.$ For real $\beta,$ we take
$U_{\beta}=\left\{  x\in A\mid\beta<x^{\ast}\left(  x\right)  \right\}  $ and
we have that
\[
A\subset%
%TCIMACRO{\dbigcup \limits_{\beta>\alpha}}%
%BeginExpansion
{\displaystyle\bigcup\limits_{\beta>\alpha}}
%EndExpansion
U_{\beta}\subset%
%TCIMACRO{\dbigcup \limits_{\beta>\alpha}}%
%BeginExpansion
{\displaystyle\bigcup\limits_{\beta>\alpha}}
%EndExpansion
\left(  U_{\beta}+K\right)  .
\]
Of course, the sets $\left(  U_{\beta}\right)  $ are weakly open, so by the
weakly $K-$compactness of $A,$ there exist a $n\in\mathbb{N}\setminus\left\{
0\right\}  $ and $\left(  \beta_{i}\right)  _{i\in\overline{1,n}}%
\subset\left(  \alpha,\infty\right)  $ such that%
\[
A\subset%
%TCIMACRO{\dbigcup \limits_{i\in\overline{1,n}}}%
%BeginExpansion
{\displaystyle\bigcup\limits_{i\in\overline{1,n}}}
%EndExpansion
\left(  U_{\beta_{i}}+K\right)  .
\]
Consider $\beta=\min\left\{  \beta_{i}\mid i\in\overline{1,n}\right\}  ,$ and
we get $A\subset U_{\beta}+K.$ Therefore, for all $a\in A$ there is $y\in
U_{\beta}$ and $k\in K$ such that $a=y+k,$ whence $x^{\ast}\left(  a\right)
=x^{\ast}\left(  y\right)  +x^{\ast}\left(  k\right)  >\beta.$ Consequently,
$\alpha\geq\beta,$ and this is a contradiction.\hfill$\square$

\begin{pr}
\label{prop_canc_sol}Suppose that $A,B,C\subset X\ $are nonempty sets, $K\ $is
solid, $C$ is weakly $K-$compact, and
\[
A+C\subset C+B+\operatorname*{int}K.
\]
Then
\[
A\subset\operatorname*{conv}B+\operatorname*{int}K.
\]

\end{pr}

\noindent\textbf{Proof. }As above, it is enough to prove that $0\in
\operatorname*{conv}B+\operatorname*{int}K$ provided $C\subset
C+B+\operatorname*{int}K.$ Observe that $\operatorname*{conv}%
B+\operatorname*{int}K$ is convex and open. Suppose, by way of contradiction,
that $0\notin\operatorname*{conv}B+\operatorname*{int}K.$ Then, using a
well-known separation result, there is $x^{\ast}\in X^{\ast}\setminus\left\{
0\right\}  $ such that
\[
0<x^{\ast}\left(  b+k\right)  ,\text{ }\forall b\in\operatorname*{conv}%
B,\text{ }\forall k\in\operatorname*{int}K.
\]
In particular, it follows that $x^{\ast}\in K^{+}.$ From Lemma \ref{lm_min},
$x^{\ast}$ attains its minimum on $C$ at a point denoted $\overline{c}$, that
is
\[
x^{\ast}\left(  \overline{c}\right)  \leq x^{\ast}\left(  c\right)  ,\text{
}\forall c\in C.
\]
Then there are $c\in C,$ $b\in B$ and $k\in\operatorname*{int}K$ such that
$\overline{c}=c+b+k,$ so%
\[
x^{\ast}\left(  \overline{c}\right)  =x^{\ast}\left(  c\right)  +x^{\ast
}\left(  b+k\right)  \geq x^{\ast}\left(  \overline{c}\right)  +x^{\ast
}\left(  b+k\right)  >x^{\ast}\left(  \overline{c}\right)  ,
\]
and this is a contradiction.\hfill$\square$

\begin{rmk}
\label{rmk_prop}Observe that if $B$ is $K-$convex, then $\operatorname*{conv}%
B+\operatorname*{int}K=B+\operatorname*{int}K.$ On the other hand, the
implication
\[
\left(  A\subset B+\operatorname*{int}K\right)  \implies\left(  A+C\subset
C+B+\operatorname*{int}K\right)
\]
always holds.
\end{rmk}

If $K$ is not necessarily solid, we have the following similar result.

\begin{pr}
\label{prop_comp_open_K}Suppose that $A,B,C\subset X\ $are nonempty sets, $C$
is weakly $K-$compact, $B$ is open, and
\[
A+C\subset C+B+K.
\]
Then
\[
A\subset\operatorname*{conv}B+K.
\]

\end{pr}

\noindent\textbf{Proof. }Again, it is enough to prove that $0\in
\operatorname*{conv}B+K$ provided $C\subset C+B+K.$ Observe that
$\operatorname*{conv}B+K=\operatorname*{conv}\left(  B+K\right)  .$ Moreover,
$B+K$ is open and the proof is as above.\hfill$\square$

\bigskip

An easy consequence that can have interesting applications to the study of set
optimization problems with particular data is recorded next. It basically says
that under certain topological assumptions some set-orders (see \cite{KTY},
for instance) are insensitive to the additions with weakly $K-$compact sets.

\begin{cor}
Suppose that $A,B,C\subset X\ $are nonempty sets, $C$ is weakly $K-$compact,
$B$ is open and convex, and
\[
A\not \subset B+K.
\]
Then
\[
A+C\not \subset C+B+K.
\]

\end{cor}

For nonempty subsets $A,B$ of $X$ one defines (see \cite{GU}, for instance)%
\[
A\overset{\bullet}{-}B=\left\{  x\in X\mid x+B\subset A\right\}
=\bigcap_{b\in B}\left(  A-b\right)  .
\]

\begin{pr}
\label{prop_comp_open_K_diff}Suppose that $A,B,C,D\subset X\ $are nonempty
sets, $D$ is weakly $K-$compact, $B$ is open, and
\[
A+C\subset D+B+K.
\]
Then
\[
A+C\overset{\bullet}{-}D\subset\operatorname*{conv}B+K.
\]

\end{pr}

\noindent\textbf{Proof. }One has $D+C\overset{\bullet}{-}D\subset C,$ so the
hypothesis means, in particular, that%
\[
A+D+C\overset{\bullet}{-}D\subset D+B+K
\]
and one applies Proposition \ref{prop_comp_open_K} to get the result.\hfill
$\square$

\begin{rmk}
Notice that $0\in C\overset{\bullet}{-}D\ $if and only if $D\subset C,$ so,
actually, Proposition \ref{prop_comp_open_K_diff} improves Proposition
\ref{prop_comp_open_K}.
\end{rmk}

Notice that one cannot drop $\operatorname*{conv}$ in the conclusion of the
above results. We give a very simple example to illustrate this for
Proposition \ref{prop_comp_open_K}.

\begin{examp}
Take $X=\mathbb{R}^{2},$ $K=[0,\infty)\times\left\{  0\right\}  ,$ $A=\left\{
0\right\}  ,$ $C=\left\{  \left(  1,2\right)  ,\left(  2,1\right)  \right\}  $
and $B=\left\{  \left(  0,-1\right)  ,\left(  -2,1\right)  \right\}  +B\left(
0_{\mathbb{R}^{2}},10^{-1}\right)  .$
\end{examp}

However, using a nonconvex separation functional (see \cite[Theorem
2.3.1]{GRTZ}) we can provide some results where one can simply put $B$ instead
of $\operatorname*{conv}B$ in the conclusion of a cancellation law.

In our notation, if $\operatorname*{int}K\neq\emptyset,$ one choose
$e\in\operatorname*{int}K$ and the Gerstewitz (Tammer) scalarizing functional
is $\varphi_{K,e}:X\rightarrow\mathbb{R}$
\begin{equation}
\varphi_{K,e}(x)=\inf\{t\in\mathbb{R}\mid x\in te-K\}. \label{eq. funct.}%
\end{equation}

For easy reference we denote $\varphi_{K,e}$ by $\varphi.$ We list here the
properties of $\varphi$ we use in the sequel (see \cite[Theorem 2.3.1]{GRTZ}).

\begin{thm}
\label{thm. phi} The functional $\varphi$ in (\ref{eq. funct.}) has the
following properties:

(i) $\varphi$ is continuous, sublinear, strictly$-\operatorname{int}%
K-$monotone, $K-$monotone;

(ii) for every $\lambda\in\mathbb{R}$ and $x\in X$,
\begin{equation}
\{u\in X\mid\varphi(u)\leq\lambda\}=\lambda e-K, \label{eq. level}%
\end{equation}%
\[
\{u\in X\mid\varphi(u)<\lambda\}=\lambda e-\operatorname*{int}K,
\]
and
\begin{equation}
\varphi(x+\lambda e)=\varphi(x)+\lambda. \label{eq. lin}%
\end{equation}

\end{thm}

\begin{lm}
\label{lm_min_fi}If $A$ is weakly $K-$compact, then the functional $\varphi$
achieves its minimum on $A.$
\end{lm}

\noindent\textbf{Proof. }The proof is similar with that of Lemma \ref{lm_min}.
Denote $\alpha=\inf\left\{  \varphi\left(  a\right)  \mid a\in A\right\}
\ $and suppose that $\alpha<\varphi^{\ast}\left(  a\right)  $ for all $a\in
A.$ For real $\beta,$ we take
\[
U_{\beta}=\left\{  x\in A\mid\beta<\varphi\left(  x\right)  \right\}
\]
and we observe that, actually $U_{\beta}=X\setminus\left(  \beta e-K\right)
,$ and these sets are weakly open (by the Mazur Theorem, $K$ is weakly closed)
and we have that
\[
A\subset%
%TCIMACRO{\dbigcup \limits_{\beta>\alpha}}%
%BeginExpansion
{\displaystyle\bigcup\limits_{\beta>\alpha}}
%EndExpansion
U_{\beta}\subset%
%TCIMACRO{\dbigcup \limits_{\beta>\alpha}}%
%BeginExpansion
{\displaystyle\bigcup\limits_{\beta>\alpha}}
%EndExpansion
\left(  U_{\beta}+K\right)  .
\]
Applying the weakly $K-$compactness of $A,$ as in the mentioned result, we get
a number $\beta>\alpha$ such that $A\subset U_{\beta}+K.$ Therefore, for all
$a\in A$ there is $y\in U_{\beta}$ and $k\in K$ such that $a=y+k,$ whence, by
$K-$monotonicity of $\varphi,$ we get $\varphi\left(  a\right)  =\varphi
\left(  y+k\right)  \geq\varphi\left(  y\right)  >\beta,$ which provides a
contradiction.\hfill$\square$

\bigskip

We can present now the announced cancellation law.

\begin{pr}
\label{prop_comp_open_K_nonconvex}Suppose that $A,B\subset X\ $are nonempty
sets, $K$ is solid, and $B$ is $K-$compact. Take $e\in\operatorname*{int}K.$
If
\[
A\not \subset B+K,
\]
then there is $\rho>0$ such that for every nonempty set $C$ which is weakly
$K-$compact and satisfies the inclusion $C\subset\left(  \rho e-K\right)
\cap\left(  -\rho e+K\right)  $ one has
\[
A+C\not \subset C+B+K.
\]

\end{pr}

\noindent\textbf{Proof. }As usual, consider $A=\left\{  0\right\}  ,$ so we
know that $B\cap-K=\emptyset,$ which is, of course, equivalent to $B\subset
X\setminus-K.$ According to \cite[Lemma 4.5]{DF2023-Subd}, there is
$\varepsilon>0$ such that $B-\varepsilon e\subset X\setminus-K$ (the
$K-$compactness of $B$ is required for this). This means that for all $b\in
B,$ $\varphi\left(  b-\varepsilon e\right)  =\varphi\left(  b\right)
-\varepsilon>0,$ whence $\varphi\left(  b\right)  >\varepsilon.$ Take
$\rho=2^{-1}\varepsilon$, a constant for which we show the conclusion.
Consider a nonempty set $C$ which is weakly $K-$compact and $C\subset\left(
\rho e-K\right)  \cap\left(  -\rho e+K\right)  .$ We get that $\varphi\left(
c\right)  \leq\rho$ for all $c\in C\cup\left(  -C\right)  .$ By Lemma
\ref{lm_min_fi}, there is $\overline{c}\in C,$ where $\varphi$ achieves its
minimum on $C.$ Supposing that $C\subset C+B+K$, there are $c\in C,$ $b\in B$
and $k\in K$ such that $\overline{c}=c+b+k.$ Therefore, using the properties
of $\varphi,$%
\[
\varphi\left(  \overline{c}\right)  =\varphi\left(  c+b+k\right)  \geq
\varphi\left(  c+b\right)  \geq\varphi\left(  b\right)  -\varphi\left(
-c\right)  >\varepsilon-\varphi\left(  -c\right)  \geq\rho\geq\varphi\left(
c\right)  \geq\varphi\left(  \overline{c}\right)  .
\]
Obviously, this is a contradiction, so the conclusion takes place.\hfill
$\square$

\bigskip

We deal now with a cancellation law that works on finite dimensional spaces
and is inspired by \cite[Proposition 5.2]{GKKU}.

\begin{pr}
\label{pr_Rad_finit}Suppose that $X$ is finite dimensional, $A,B,C\subset
X\ $are nonempty sets, $C$ is $K-$sequentially compact, and
\[
A+C\subset C+B+K.
\]
Then
\[
A\subset\operatorname*{conv}B+K.
\]

\end{pr}

\noindent\textbf{Proof. }For $K=\left\{  0\right\}  $ the result is
\cite[Proposition 5.2]{GKKU}, whence we consider that $K$ is proper. (Notice
that the original proof of the mentioned result is based on a separation
theorem with respect to the lexicographical order from \cite{ML}.) Again, it
is enough to prove for $A=\left\{  0\right\}  $. Also by using a translation
argument, one can assume that $0\in C\ $and, moreover, one can consider $B$ to
be convex.

So, we have to show that $0\in B+K$ provided $0\in C\subset C+B+K,$ $B$ is
convex and $C$ is $K-$sequentially compact.

We proceed by induction. Consider the proposition:

$P\left(  n\right)  :$ if $X$ has dimension $n,$ for all $B,C,K\subset X$ with
$B,C\neq\emptyset,$ $K$ a pointed proper convex and closed cone such that $C$
is $K-$sequentially compact, $B$ is convex and $0\in C\subset C+B+K,$ one has
$0\in B+K$.

Take $n=1.$ Then, $X=\mathbb{R}x$ with $x\in S_{X}.$ Without loss of
generality, we take $K=[0,\infty)x.$ Since $C$ is $K-$sequentially compact, it
is $K-$bounded and $K-$closed (see \cite{DF2022}), so $\inf\left\{  c\mid
cx\in C\right\}  =\min\left\{  c\mid cx\in C\right\}  \in\mathbb{R}.$ Denoting
$\inf\left\{  c\mid cx\in C\right\}  $ by $\overline{c}$ and using the
assumption, there are $c,b,k\in\mathbb{R}$ with $cx\in C,$ $bx\in B,$ $kx\in
K$ such that
\[
\overline{c}=c+b+k\geq\overline{c}+b+k,
\]
which means that $b\leq0,$ whence $0\in B+[0,\infty)x=B+K.$

Suppose now that $P\left(  n\right)  $ is true and prove that $P\left(
n+1\right)  $ is true. Suppose that $0\notin B+K.$ By a standard separation
argument, there is $x\in X\setminus\left\{  0\right\}  $ such that%
\[
0\leq\left\langle x,b+k\right\rangle ,\text{ }\forall b\in B,\text{ }\forall
k\in K,
\]
where $\left\langle \cdot,\cdot\right\rangle $ stands for the usual inner
product of $X.$ Surely, this implies as well that
\[
0\leq\left\langle x,b\right\rangle ,\text{ }\forall b\in B\text{ and }%
0\leq\left\langle x,k\right\rangle ,\text{ }\forall k\in K.
\]
Now, the $K-$sequentially compactness of $C$ and the latter property of $x$
ensure (by Lemma \ref{lm_min}, for instance) that there is $\overline{c}\in C$
such that%
\[
\left\langle x,\overline{c}\right\rangle \leq\left\langle x,c\right\rangle
,\text{ }\forall c\in C.
\]
Now, there are $c\in C,$ $b\in B$ and $k\in K$ such that $\overline{c}=c+b+k,$
so $0=c-\overline{c}+b+k\in C-\overline{c}+B+K.$ Denote $C^{\prime
}=C-\overline{c}$ and observe that $\left\langle x,c^{\prime}\right\rangle
\geq0$ for all $c^{\prime}\in C^{\prime}.$ We have
\[
0=\left\langle x,c-\overline{c}+b+k\right\rangle =\left\langle x,c-\overline
{c}\right\rangle +\left\langle x,b\right\rangle +\left\langle x,k\right\rangle
.
\]
We get from here that $\left\langle x,c-\overline{c}\right\rangle
=\left\langle x,b\right\rangle =\left\langle x,k\right\rangle =0.$ Consider
now the $n-$dimensional subspace of $X$ as $Z=\left\{  z\in X\mid\left\langle
x,z\right\rangle =0\right\}  $ and denote $C_{1}=C^{\prime}\cap Z,$
$B_{1}=B\cap Z$ and $K_{1}=K\cap Z.$ Therefore, $0=\overline{c}-\overline
{c}\in C_{1},$ $b\in B_{1},$ $k\in K_{1}.$ Take $y\in C_{1}.$ There are
$c^{\prime}\in C^{\prime},$ $b^{\prime}\in B,$ $k^{\prime}\in K$ such that
$y=c^{\prime}+b^{\prime}+k^{\prime}.$ We have
\[
0=\left\langle x,c^{\prime}+b^{\prime}+k^{\prime}\right\rangle =\left\langle
x,c^{\prime}\right\rangle +\left\langle x,b^{\prime}\right\rangle
+\left\langle x,k^{\prime}\right\rangle .
\]
We deduce that $c^{\prime}\in C_{1},$ $b^{\prime}\in B_{1},$ $k^{\prime}\in
K_{1},$ so $y\in C_{1}+B_{1}+K_{1}.$ We get $0\in C_{1}\subset C_{1}%
+B_{1}+K_{1}.$ Clearly, $B_{1}$ is a convex set. In order to apply $P\left(
n\right)  $ we have to prove that $C_{1}$ is $K_{1}-$sequentially compact.
Take $\left(  c_{n}\right)  \subset C_{1}=C^{\prime}\cap Z.$ In particular,
$\left(  c_{n}\right)  \subset C^{\prime}$ and $C^{\prime}$ is $K-$%
sequentially compact, being a translation of a $K-$sequentially compact set.
Therefore, there is $\left(  k_{n}\right)  \subset K$ such that, on a
subsequence, $c_{n}-k_{n}\rightarrow u\in C^{\prime}$. We have that
\[
\left\langle x,c_{n}-k_{n}\right\rangle =\left\langle x,c_{n}\right\rangle
-\left\langle x,k_{n}\right\rangle \rightarrow\left\langle x,u\right\rangle .
\]
But $\left\langle x,c_{n}\right\rangle =0$ and $\left\langle x,k_{n}%
\right\rangle \geq0$ for all $n$, while $\left\langle x,u\right\rangle \geq0.$
We deduce that $\left\langle x,u\right\rangle =0$, whence $u\in C_{1}$ and
$\left\langle x,k_{n}\right\rangle \rightarrow0$ (on a subsequence). The
latter relation shows that $d\left(  k_{n},Z\right)  \rightarrow0$ and using
the metric regularity properties of sets proved in \cite[Theorems 3.9,
3.17]{BB} we get that $d\left(  k_{n},K_{1}\right)  \rightarrow0.$ For all $n$
denote by $k_{n}^{\prime}$ the projection of $k_{n}$ on $K_{1}$ and we have
that $k_{n}-k_{n}^{\prime}\rightarrow0.$ Therefore $c_{n}-k_{n}^{\prime
}\rightarrow u\in C_{1}$ and we conclude that $C_{1}$ is $K_{1}-$sequentially
compact. By $P\left(  n\right)  ,$ we deduce that $0\in B_{1}+K_{1}\subset
B+K.$ This is a contradiction. The conclusion follows.\hfill$\square$

\section{Some applications}

In this short section we derive some consequences of the conic cancellation
laws of the previous section. The first results are in the spirit of
\cite[Lemma 3]{Rad} (see also \cite{Bat}) and, first of all, some remarks are
in order.

\begin{rmk}
\label{rmk_bounded_K}If $A\subset X$ is a $K-$bounded set and $B\subset X$ is
nonempty, then $e\left(  A+K,B+K\right)  \in\mathbb{R}.$ Indeed, denoting by
$M$ a bounded set for which $A\subset M+K,$ and taking $b\in B,$ we have
\begin{align*}
e\left(  A+K,B+K\right)   &  =\sup\left\{  d\left(  a+c,B+K\right)  \mid a\in
A,c\in K\right\} \\
&  \leq\sup\left\{  d\left(  a+c,B+K\right)  \mid a\in M,c\in K\right\} \\
&  \leq\sup\left\{  d\left(  a+c,b+K\right)  \mid a\in M,c\in K\right\} \\
&  \leq\sup\left\{  \left\Vert a-b\right\Vert \mid a\in M\right\}  <+\infty.
\end{align*}
Moreover, using relations (\ref{rel1}) and (\ref{rel2}), it is not difficult
to see that%
\[
e(A,B+K)=e(A+K,B+K)=e\left(  A,\operatorname{cl}\left(  B+K\right)  \right)
=e\left(  A,\operatorname{cl}B+K\right)  =e\left(  \operatorname{cl}%
A+K,\operatorname{cl}B+K\right)  .
\]

\end{rmk}

A conic counterpart of \cite[Lemma 3]{Rad} reads as follows.

\begin{pr}
\label{pr_Rlema3_1}Suppose that $A,B,C\subset X\ $are nonempty sets, $C$ is
weakly $K-$compact, and $\operatorname*{cl}B$ is $K-$convex. Then%
\[
e\left(  A,B+K\right)  =e\left(  A+C,B+K+C\right)  .
\]

\end{pr}

\noindent\textbf{Proof. }As mentioned,%
\[
e\left(  A,B+K\right)  =\inf\left\{  \alpha>0\mid A\subset B+\alpha
B_{X}+K\right\}  ,
\]
while%
\[
e\left(  A+C,B+K+C\right)  =\inf\left\{  \alpha>0\mid A+C\subset B+\alpha
B_{X}+K+C\right\}  .
\]
But, for all $\alpha>0,$ the set $B+\alpha B_{X}+K$ is open and convex, the
latter assertion being based on the equality%
\[
B+\alpha B_{X}+K=\operatorname*{cl}B+\alpha B_{X}+K,\text{ }\forall\alpha>0,
\]
which is a consequence of $K-$convexity of $\operatorname*{cl}B.$ So according
to Proposition \ref{prop_comp_open_K}, the relations $A\subset B+\alpha
B_{X}+K$ and $A+C\subset B+\alpha B_{X}+K+C$ are equivalent. The conclusion
ensues.\hfill$\square$

\begin{rmk}
Actually, it is not difficult to observe that the convexity of
$\operatorname*{cl}B$ is equivalent to the convexity of $B+\alpha B_{X}$ for
any $\alpha>0.$
\end{rmk}

\bigskip

We record the following consequence of the Propositions \ref{prop_canc_sol}
and \ref{pr_Rlema3_1}.

\begin{cor}
Suppose that $K\ $is solid, $A,B,C\subset X\ $are nonempty, weakly
$K-$compact, and $K-$convex sets. Then:

(i) relation $A+C+\operatorname*{int}K=C+B+\operatorname*{int}K$ implies
$A+\operatorname*{int}K=B+\operatorname*{int}K;$

(ii) $e\left(  A+\operatorname*{int}K,B+\operatorname*{int}K\right)  =e\left(
A+C+\operatorname*{int}K,B+C+\operatorname*{int}K\right)  .$
\end{cor}

\noindent\textbf{Proof. }(i) Taking into account that $\operatorname*{int}%
K=\operatorname*{int}K+K$ this item is a consequence of Proposition
\ref{prop_canc_sol}.

(ii) Since $A$ is weakly $K-$compact then $A+K$ is weakly closed (see
\cite[Proposition 3.3, p. 14]{Luc}), hence strongly closed and the same can be
said about $B$ and $C.$ Moreover, it is easy to see that $e\left(
A,B+\operatorname*{int}K\right)  =e\left(  A,B+K\right)  $ and then we can
apply Proposition \ref{pr_Rlema3_1}.\hfill$\square$

\begin{rmk}
This corollary ensures the essential properties a semigroup should have in
order to be embedded into a vector space (according to \cite[Theorem 1]{Rad})
for the semigroup of the sets of the form $A+\operatorname*{int}K,$ where $A$
is nonempty, weakly $K-$compact, and $K-$convex.
\end{rmk}

\bigskip

On the basis of Proposition \ref{pr_Rad_conic} one can get a result similar to
Proposition \ref{pr_Rlema3_1}.

\begin{pr}
Suppose that $A,B,C\subset X\ $are nonempty sets such that $C$ is $K-$bounded
and $B+\alpha D_{X}$ is $K-$convex and $K-$closed for any $\alpha>0.$ Then%
\[
e\left(  A,B+K\right)  =e\left(  A+C,C+B+K\right)  .
\]

\end{pr}

\noindent\textbf{Proof. }We use
\[
e\left(  A,B+K\right)  =\inf\left\{  \alpha>0\mid A\subset B+\alpha
D_{X}+K\right\}  ,
\]
and%
\[
e\left(  A+C,B+K+C\right)  =\inf\left\{  \alpha>0\mid A+C\subset B+\alpha
D_{X}+K+C\right\}
\]
and Proposition \ref{pr_Rad_conic} to conclude that under our assumptions, the
relations $A\subset B+\alpha D_{X}+K$ and $A+C\subset B+\alpha D_{X}+K+C$ are
equivalent.\hfill$\square$

\begin{rmk}
Observe that if $B$ is $K-$convex and $K-$sequentially compact, and $D_{X}$ is
$K-$closed, then $B+\alpha D_{X}$ is $K-$convex and $K-$closed for any
$\alpha>0.$
\end{rmk}

\begin{cor}
Suppose that $A,B,C\subset X\ $are nonempty, $K-$sequentially compact, and
$K-$convex sets. If $D_{X}$ is $K-$closed, then:

(i) relation $A+C+K=C+B+K$ implies $A+K=B+K;$

(ii) $e\left(  A+K,B+K\right)  =e\left(  A+C+K,B+C+K\right)  .$
\end{cor}

\bigskip

Next, we consider another normed vector space $Z$ and $F:Z\rightrightarrows X$
a set-valued map. One considers as well the associated epigraphical set-valued
map $\operatorname*{Epi}F:Z\rightrightarrows X$ given by $\operatorname*{Epi}%
F\left(  z\right)  =F\left(  z\right)  +K.$ Recall the following notion from
\cite{DS23-JOGO}.

\begin{df}
\label{subdiff}The Fr\'{e}chet subdifferential of $F$ at $\overline{z}\in Z$
is%
\begin{equation}
\widehat{\partial}F\left(  \overline{z}\right)  =\left\{  T\in B\left(
Z,X\right)  \mid\lim_{z\rightarrow\overline{z}}\frac{e\left(
\operatorname*{Epi}F\left(  z\right)  ,\operatorname*{Epi}F\left(
\overline{z}\right)  +T\left(  z-\overline{z}\right)  \right)  }{\left\Vert
z-\overline{z}\right\Vert }=0\right\}  . \label{Fr}%
\end{equation}
Equivalently, $T\in\widehat{\partial}F\left(  \overline{z}\right)  $ iff $T\in
B\left(  Z,X\right)  $ and%
\begin{equation}
\forall\varepsilon>0,\exists\delta>0,\forall z\in B\left(  \overline{z}%
,\delta\right)  :\operatorname*{Epi}F\left(  z\right)  \subset
\operatorname*{Epi}F\left(  \overline{z}\right)  +T\left(  z-\overline
{z}\right)  +\varepsilon\left\Vert z-\overline{z}\right\Vert D_{X}.
\label{Fr2}%
\end{equation}

Similarly, we can define the upper subdifferential of $F$ at $\overline{z}$ as
follows%
\begin{equation}
\widehat{\partial}^{+}F\left(  \overline{z}\right)  =\left\{  T\in B\left(
Z,X\right)  \mid\lim_{z\rightarrow\overline{z}}\frac{e\left(
\operatorname*{Epi}F\left(  \overline{z}\right)  +T\left(  z-\overline
{z}\right)  ,\operatorname*{Epi}F\left(  z\right)  \right)  }{\left\Vert
z-\overline{z}\right\Vert }=0\right\}  . \label{uFr}%
\end{equation}

\end{df}

Remark that, if $f:X\rightarrow\mathbb{R}$ is a function, then relation
(\ref{Fr2}) can be equivalently written as%
\[
\forall\varepsilon>0,\exists\delta>0,\forall z\in B\left(  \overline{z}%
,\delta\right)  :f\left(  z\right)  -f\left(  \overline{z}\right)  -z^{\ast
}\left(  z-\overline{z}\right)  \geq-\varepsilon\left\Vert z-\overline
{z}\right\Vert ,
\]
i.e., $z^{\ast}\in\widehat{\partial}f\left(  \overline{z}\right)  ,$ where
$\widehat{\partial}f\left(  \overline{z}\right)  $ denotes the usual
Fr\'{e}chet subdifferential of $f$ at $\overline{z}$ (see \cite{Morduk2006}).
Comments of the same kind are in order for $\widehat{\partial}^{+}.$

We derive some new calculus rules for this Fr\'{e}chet subdifferential, while
for upper subdifferential similar results holds.

\begin{pr}
If $A\subset X$ is a nonempty weakly $K-$compact set, $\overline{z}\in X$ and
$\operatorname*{cl}F\left(  \overline{z}\right)  $ is $K-$convex, then
$\widehat{\partial}\left(  F\left(  \cdot\right)  +A\right)  \left(
\overline{z}\right)  =\widehat{\partial}F\left(  \overline{z}\right)  .$
\end{pr}

\noindent\textbf{Proof. }Using Proposition \ref{pr_Rlema3_1}, one has, in the
notation of Definition \ref{subdiff}, that%
\[
e\left(  \operatorname*{Epi}F\left(  z\right)  ,\operatorname*{Epi}F\left(
\overline{z}\right)  +T\left(  z-\overline{z}\right)  \right)  =e\left(
\operatorname*{Epi}F\left(  z\right)  +A,\operatorname*{Epi}F\left(
\overline{z}\right)  +A+T\left(  z-\overline{z}\right)  \right)  ,
\]
so the required equality holds$.$\hfill$\square$

\begin{df}
One says that $F$ is $K-$Lipschitz around $\overline{z}\in Z$ if there are a
neighborhood $U$ of $\overline{z}$, a constant $\ell>0$ and an element $e\in
S_{Z}\cap K$ such that for every $z^{\prime},z^{\prime\prime}\in U,$%
\[
F(z^{\prime\prime})+\ell\left\Vert z^{\prime\prime}-z^{\prime}\right\Vert
e\subset F(z^{\prime})+K.
\]

\end{df}

Next, we additionally use the notion of normal cone (see \cite[Definition
2.1.21]{GRTZ}) by means of some of its characterizations (see \cite[Theorem
2.2.10]{GRTZ}). Rather than the formal definition we present a
characterization of this concept we use in the sequel.

\begin{pr}
\label{prop_con_nor}The cone $K$ is normal iff there exists $\alpha>0$ such
that $\left\Vert z\right\Vert \leq\alpha\left\Vert y\right\Vert $ whenever
$z,y\in K,$ $y-z\in K$.
\end{pr}

\begin{pr}
If $F$ is $K-$Lipschitz (with constant $\ell$ and element $e\in S_{Z}\cap K$)
around $\overline{z}\in X,$ has weakly $K-$compact values, and $K$ is normal,
then there is $\alpha>0$ such that for all $T\in\widehat{\partial}F\left(
\overline{z}\right)  $ and $u\in S_{Z}\cap T^{-1}\left(  K\cup\left(
-K\right)  \right)  ,$ one has
\[
\left\Vert Tu\right\Vert \leq\alpha\ell\left\Vert e\right\Vert .
\]

\end{pr}

\noindent\textbf{Proof. }Take $T\in\widehat{\partial}F\left(  \overline
{z}\right)  .$ Taking into account the Lipschitz property, for all
$\varepsilon>0$ there is $\delta>0$ such that for all $z\in B\left(
\overline{z},\delta\right)  ,$%
\begin{align*}
F\left(  z\right)   &  \subset F\left(  \overline{z}\right)  +T\left(
z-\overline{z}\right)  +\varepsilon\left\Vert z-\overline{z}\right\Vert
B_{X}+K\\
&  \subset F\left(  z\right)  -\ell\left\Vert z-\overline{z}\right\Vert
e+T\left(  z-\overline{z}\right)  +\varepsilon\left\Vert z-\overline
{z}\right\Vert B_{X}+K.
\end{align*}
Using Proposition \ref{prop_comp_open_K}, one has for all $z\in B\left(
\overline{z},\delta\right)  ,$%
\[
0\in T\left(  z-\overline{z}\right)  +\left\Vert z-\overline{z}\right\Vert
\left(  -\ell e+\varepsilon B_{X}\right)  +K,
\]
which means that for all $z\in S_{Z},$%
\[
T\left(  z\right)  \in-\ell e+\varepsilon B_{X}+K.
\]
Fix $u\in S_{Z}\cap T^{-1}\left(  K\cup\left(  -K\right)  \right)  $ and
observe that this set is symmetric. Taking $u$ or $-u$ we can suppose, without
loss of generality, that $T\left(  u\right)  \in K$ and
\[
\ell e-T\left(  u\right)  \in\varepsilon B_{X}+K.
\]
Notice that this holds for all $\varepsilon>0.$ Therefore, $\ell e-T\left(
u\right)  \in\operatorname{cl}K=K,$ and the normality of $K$ implies that
$\left\Vert Tu\right\Vert \leq\alpha\ell\left\Vert e\right\Vert ,$ where
$\alpha$ is the constant from Proposition \ref{prop_con_nor}.\hfill$\square$

\begin{rmk}
Notice that in \cite{DF2023-Subd}, besides the study of the Fr\'{e}chet
subdifferentials, a limiting type (that is, Mordukhovich type: see
\cite{Morduk2006}) subdifferential was introduced using, as usual, a limiting
procedure. It is easy to see that the above results can be readily extended to
this type of generalized subgradients.
\end{rmk}

\bigskip

Finally, we deal with a sharp type solution for a set optimization problem.
Let $M\subset Z$ be a closed set. Recall that the Bouligand tangent cone to
$M$ at $\overline{z}\in M$ is the set%
\[
T_{B}\left(  M,\overline{z}\right)  :=\left\{  u\in X\mid\exists\left(
t_{n}\right)  \downarrow0,\text{ }\exists\left(  u_{n}\right)  \rightarrow
u,\text{ }\forall n\in\mathbb{N},\text{ }\overline{z}+t_{n}u_{n}\in M\right\}
.
\]

\begin{df}
In the above notation if $K$ is solid, we say that $\overline{z}\in M$ is a
sharp weak minimum for $F$ on $M$ if there is $e\in K\setminus\left\{
0\right\}  $ such that for all $z\in M,$
\[
F\left(  \overline{z}\right)  \not \subset F\left(  z\right)  -\mu\left\Vert
z-\overline{z}\right\Vert e+\operatorname*{int}K.
\]

\end{df}

This notion corresponds to a well-known concept of solution in scalar and
vectorial optimization problems (see, for instance, \cite{W}, \cite{FJ}) and
is stronger than the so-called $\ell-$minimum studied, for instance in
\cite{DS23-JOGO}. We derive a necessary optimality condition for this kind of solution.

\begin{pr}
\label{pr_sharp_d+}If $\overline{z}$ is sharp weak minimum for $F$ on $M$ then
for all $T\in\widehat{\partial}^{+}F\left(  \overline{z}\right)  $ and $u\in
T_{B}\left(  M,\overline{z}\right)  ,$%
\[
T\left(  u\right)  \notin\mu\left\Vert u\right\Vert e-\operatorname*{int}K.
\]

\end{pr}

\noindent\textbf{Proof. }Suppose the conclusion is not true, so there are
$T\in\widehat{\partial}^{+}F\left(  \overline{z}\right)  $ and $u\in
T_{B}\left(  M,\overline{z}\right)  $ such that $T\left(  u\right)  \in
\mu\left\Vert u\right\Vert e-\operatorname*{int}K.$ Then there are as well
$\left(  t_{n}\right)  \downarrow0,$ $\left(  u_{n}\right)  \rightarrow u$
such that for all $n,$ $\overline{z}+t_{n}u_{n}\in M.$ Take $\varepsilon>0$
such that $T\left(  u\right)  +4\varepsilon D_{X}\subset\mu\left\Vert
u\right\Vert e-\operatorname*{int}K.$ So, for $n$ large enough, we have from
the definitions of $\widehat{\partial}^{+}F\left(  \overline{z}\right)  $ and
sharp weak minimum%
\begin{align}
F\left(  \overline{z}\right)   &  \subset F\left(  \overline{z}+t_{n}%
u_{n}\right)  -t_{n}T\left(  u_{n}\right)  +\varepsilon t_{n}D_{X}%
+K\label{rel_con}\\
F\left(  \overline{z}\right)   &  \not \subset F\left(  \overline{z}%
+t_{n}u_{n}\right)  -t_{n}\mu\left\Vert u_{n}\right\Vert e+\operatorname*{int}%
K.\nonumber
\end{align}
But, again, for large $n,$ since $T$ is continuous and $u_{n}\rightarrow u$ we
have
\[
-T\left(  u_{n}\right)  +2\varepsilon D_{X}\subset-T\left(  u\right)
+3\varepsilon D_{X}%
\]
and
\[
\mu\left\Vert u\right\Vert e-\operatorname*{int}K\subset\mu\left\Vert
u_{n}\right\Vert e+\varepsilon D_{X}-\operatorname*{int}K.
\]
From the choice of $\varepsilon,$ for every large but fixed $n,$ we get%
\[
T\left(  u\right)  +4\varepsilon D_{X}=T\left(  u\right)  +3\varepsilon
D_{X}+\varepsilon D_{X}\subset\mu\left\Vert u\right\Vert e-\operatorname*{int}%
K\subset\mu\left\Vert u_{n}\right\Vert e+\varepsilon D_{X}-\operatorname*{int}%
K,
\]
and from R\aa dstr\"{o}m cancellation law, we deduce that%
\[
T\left(  u\right)  +3\varepsilon D_{X}\subset\operatorname{cl}\left(
\mu\left\Vert u_{n}\right\Vert e-\operatorname*{int}K\right)  =\mu\left\Vert
u_{n}\right\Vert e-K.
\]

Consequently,
\[
-T\left(  u_{n}\right)  +\varepsilon D_{X}+\varepsilon D_{X}=-T\left(
u_{n}\right)  +2\varepsilon D_{X}\subset-T\left(  u\right)  +3\varepsilon
D_{X}\subset-\mu\left\Vert u_{n}\right\Vert e+K,
\]
and, therefore,%
\[
-T\left(  u_{n}\right)  +\varepsilon D_{X}\subset-\mu\left\Vert u_{n}%
\right\Vert e+\operatorname*{int}K.
\]
This contradicts the group of relations (\ref{rel_con}). So, the conclusion
holds.\hfill$\square$

\begin{rmk}
Notice that the above result works as well for easy-to-define local
counterpart of the sharp weak minimality.
\end{rmk}

\begin{rmk}
In the context of a real-valued function $f$, when $\overline{z}$ is an
unconstrained minimum for $f$, then $0\in\widehat{\partial}f\left(
\overline{z}\right)  \neq\emptyset,$ whence for $F=f,$ $\widehat{\partial}%
^{+}F\left(  \overline{z}\right)  =\widehat{\partial}^{+}f\left(  \overline
{z}\right)  \neq\emptyset$ iff $f$ is Fr\'{e}chet differentiable, in which
case $\widehat{\partial}f\left(  \overline{z}\right)  =\widehat{\partial}%
^{+}f\left(  \overline{z}\right)  =\left\{  \nabla f\left(  \overline
{z}\right)  \right\}  $ (see \cite[Proposition 1.87]{Morduk2006}). For
instance, if we take the example of an interval set-valued map $F\left(
z\right)  =\left[  f\left(  z\right)  ,g\left(  z\right)  \right]  $ with
$z\in Z$ and $f,g:Z\rightarrow\mathbb{R}$, $f\leq g,$ in the case $f$ is
Fr\'{e}chet differentiable, then $\left\{  \nabla f\left(  \overline
{z}\right)  \right\}  =\widehat{\partial}F\left(  \overline{z}\right)
\cap\widehat{\partial}^{+}F\left(  \overline{z}\right)  .$ However, in the
setting of Proposition \ref{pr_sharp_d+}, one does not necessarily have
$0\in\widehat{\partial}F\left(  \overline{z}\right)  ,$ so $\widehat{\partial
}^{+}F\left(  \overline{z}\right)  \neq\emptyset$ is not so restrictive.

Another useful situation (see \cite{HL}, \cite{EGR}) is as follows: consider
$f_{i},g_{i}:Z\rightarrow\mathbb{R}$, $f_{i}\leq g_{i}$ with $i\in\left\{
1,2\right\}  $ and $F:Z\rightrightarrows\mathbb{R}^{2}$ given by%
\[
F\left(  z\right)  =\left[  f_{1}\left(  z\right)  ,g_{1}\left(  z\right)
\right]  \times\left[  f_{2}\left(  z\right)  ,g_{2}\left(  z\right)  \right]
.
\]
Take $K=\mathbb{R}_{+}^{2}.$ Then for all $\overline{z}\in Z,$
\[
\widehat{\partial}F\left(  z\right)  =\left\{  \left(  z_{1}^{\ast}%
,z_{2}^{\ast}\right)  \mid z_{i}^{\ast}\in\widehat{\partial}f_{i}\left(
\overline{z}\right)  ,\text{ }i\in\left\{  1,2\right\}  \right\}  ,
\]
and similarly for $\widehat{\partial}^{+}F\left(  \overline{z}\right)  .$ As
mentioned in \cite{DF2023-Subd}, the $\ell-$minimality of $\overline{z}$ for
$F$ does not force in this case $\widehat{\partial}F\left(  \overline
{z}\right)  \neq\emptyset.$
\end{rmk}

\bigskip

We end with a proposition inspired by the perturbation result for scalar
optimization problems presented in \cite[Lemma 2.1]{Sha}. The possibility to
deal with sets instead of points in this kind of assertions is driven by
cancellation laws.

\begin{pr}
Take $F,H:Z\rightrightarrows X$ be set-valued maps, $M\subset Z$ be a nonempty
closed set, and $e\in\operatorname*{int}K,$ $\mu>0,$ $L\in\left(
0,\mu\right)  ,$ $\varepsilon>0.$ Suppose that exist $\overline{z}%
,z_{\varepsilon}\in M$ such that for all $z\in M,$

(i) $F\left(  \overline{z}\right)  \not \subset F\left(  z\right)
-\mu\left\Vert z-\overline{z}\right\Vert e+\operatorname*{int}K;$

(ii) $H\left(  z\right)  +L\left\Vert z-\overline{z}\right\Vert e\subset
H\left(  \overline{z}\right)  +K;$

(iii) $\left(  F+H\right)  \left(  z\right)  \subset\left(  F+H\right)
\left(  z_{\varepsilon}\right)  -\varepsilon e+K.$

Moreover, suppose that $F\left(  z_{\varepsilon}\right)  $ is $K-$convex and
$H\left(  \overline{z}\right)  $ is weakly $K-$compact.

Then
\[
\left\Vert z_{\varepsilon}-\overline{z}\right\Vert \leq\frac{\varepsilon}%
{\mu-L}.
\]

\end{pr}

\noindent\textbf{Proof. }We know, from $(i)$ that
\[
F\left(  \overline{z}\right)  \not \subset F\left(  z_{\varepsilon}\right)
-\mu\left\Vert z_{\varepsilon}-\overline{z}\right\Vert e+\operatorname*{int}K.
\]
According to Proposition \ref{prop_canc_sol}, this is equivalent to
\[
F\left(  \overline{z}\right)  +H\left(  \overline{z}\right)  \not \subset
F\left(  z_{\varepsilon}\right)  +H\left(  \overline{z}\right)  -\mu\left\Vert
z_{\varepsilon}-\overline{z}\right\Vert e+\operatorname*{int}K.
\]
On the other hand, by $(iii)$,%
\[
F\left(  \overline{z}\right)  +H\left(  \overline{z}\right)  =\left(
F+H\right)  \left(  \overline{z}\right)  \subset\left(  F+H\right)  \left(
z_{\varepsilon}\right)  -\varepsilon e+K,
\]
whence%
\[
\left(  F+H\right)  \left(  z_{\varepsilon}\right)  -\varepsilon
e+K\not \subset F\left(  z_{\varepsilon}\right)  +H\left(  \overline
{z}\right)  -\mu\left\Vert z_{\varepsilon}-\overline{z}\right\Vert
e+\operatorname*{int}K.
\]
We deduce that (see Remark \ref{rmk_prop}, as well)%
\[
H\left(  z_{\varepsilon}\right)  -\varepsilon e+K\not \subset H\left(
\overline{z}\right)  -\mu\left\Vert z_{\varepsilon}-\overline{z}\right\Vert
e+\operatorname*{int}K.
\]
Then $(ii)$ gives%
\[
H\left(  z_{\varepsilon}\right)  +L\left\Vert z_{\varepsilon}-\overline
{z}\right\Vert e\subset H\left(  \overline{z}\right)  +K,
\]
so%
\[
H\left(  \overline{z}\right)  -L\left\Vert z_{\varepsilon}-\overline
{z}\right\Vert e-\varepsilon e+K\not \subset H\left(  \overline{z}\right)
-\mu\left\Vert z_{\varepsilon}-\overline{z}\right\Vert e+\operatorname*{int}K.
\]
This gives%
\[
\left(  \mu-L\right)  \left\Vert z_{\varepsilon}-\overline{z}\right\Vert
e-\varepsilon e\notin\operatorname*{int}K,
\]
which is
\[
\left(  \mu-L\right)  \left\Vert z_{\varepsilon}-\overline{z}\right\Vert
\leq\varepsilon,
\]
and the conclusion follows.\hfill$\square$

\bigskip

\noindent\textbf{Acknowledgements. }The authors are grateful to Professor
Constantin Z\u{a}linescu for providing an alternative proof of
\cite[Proposition 5.2]{GKKU} by mathematical induction on which the proof of
Proposition \ref{pr_Rad_finit} is based.

\noindent\textbf{Funding.} This work was supported by a grant of the Ministry
of Research, Innovation and Digitization, CNCS - UEFISCDI, project number
PN-III-P4-PCE-2021-0690, within PNCDI III.

\noindent\textbf{Data availability.} This manuscript has no associated data.

\noindent\textbf{Disclosure statement}. No potential conflict of interest was
reported by the authors.


\bigskip

\begin{thebibliography}{99}                                                                                               %


\bibitem {Bat}R.G. Batson, \emph{Extensions of R\aa dstr\"{o}m's lemma with
application to stability theory of mathematical programming}, Journal of
Mathematical Analysis and Applications, 117 (1986), 441--448.

\bibitem {BB}H.H. Bauschke, J.M. Borwein, \emph{On the convergence of von
Neumann's alternating projection algorithm for two sets}, Set-Valued Analysis,
1 (1993), 185--212.

\bibitem {BT}J. Bielawski, J. Tabor, \emph{An embedding theorem for unbounded
convex sets in a Banach space}, Demonstratio Mathematica, 42 (2009), 703--709.

\bibitem {DF2022}M. Durea, E.-A. Florea, \emph{Cone-compactness of a set and
applications to set-equilibrium problems}, submitted.

\bibitem {DF2023-Subd}M. Durea, E.-A. Florea, \emph{Subdifferential calculus
and ideal solutions for set optimization problems}\textbf{,} submitted.

\bibitem {DS23-JOGO}M. Durea, R. Strugariu, \emph{Directional derivatives and
subdifferentials for set-valued maps applied to set optimization}, Journal of
Global Optimization, 85 (2023), 687--707.

\bibitem {EGR}G. Eichfelder, T. Gerlach, S. Rockt\"{a}schel, \emph{Convexity
and continuity of specific set-valued maps and their extremal value
functions}, https://optimization-online.org/wp-content/uploads/2022/04/8883.pdf.

\bibitem {FJ}F. Flores-Baz\'{a}n, B. Jim\'{e}nez, \emph{Strict efficiency in
set-valued optimization}, SIAM Journal on Control and Optimization, 48 (2009), 881--908.

\bibitem {GKKU}J. Grzybowski, M. K\"{u}\c{c}\"{u}k, Y. K\"{u}\c{c}\"{u}k, R.
Urb\'{a}nski, \emph{Minkowski--R\aa dstr\"{o}m--H\"{o}rmander cone}, Pacific
Journal of Optimization, 10 (2014), 649--666.

\bibitem {GU}J. Grzybowski,\ R. Urb\'{a}nski, \emph{Order cancellation law in
the family of bounded convex sets}, Journal of Global Optimization, 77 (2020), 289--300.

\bibitem {GRTZ}A. G\"{o}pfert, H. Riahi, Chr. Tammer, C. Z\u{a}linescu,\emph{
Variational Methods in Partially Ordered Spaces}, Springer, Berlin, 2003.

\bibitem {HL}E. Hern\'{a}ndez, R. L\'{o}pez, \emph{Some useful set-valued maps
in set optimization}, Optimization, 66 (2017), 1273--1289.

\bibitem {KN}D. Kuroiwa, T. Nuriya, \emph{A generalized embedding vector space
in set optimization}, Nonlinear Analysis and Convex Analysis, 5 (2007), 297--303.

\bibitem {KTY}I. Kuwano, T. Tanaka, S. Yamada, \emph{Unified scalarization for
sets in set-valued optimization}, Nonlinear Analysis and Convex Analysis, 1685
(2010), 270--280.

\bibitem {Luc}D.T. Luc, \emph{Theory of Vector Optimization}\textit{,}
Springer, Berlin, 1989.

\bibitem {ML}J.E. Mart\'{\i}nez-Legaz, \emph{Exact quasiconvex conjugation},
Zeitschrift f\"{u}r Operations-Research, 27 (1983), 257--266.

\bibitem {Morduk2006}B.S. Mordukhovich, \emph{Variational Analysis and
Generalized Differentiation, Vol. I: Basic Theory, Vol. II: Applications},
Springer, Berlin, 2006.

\bibitem {Rad}H. R\aa dstr\"{o}m, \emph{An embedding theorem for spaces of
convex sets}, Proceedings of the American Mathematical Society, 3 (1952), 165--169.

\bibitem {Sch}K.D. Schmidt, \emph{Embedding theorems for classes of convex
sets}, Acta Applicandae Mathematicae, 5 (1986), 209--237.

\bibitem {Sha}A. Shapiro, \emph{Perturbation analysis of optimization problems
in Banach spaces}, Numerical Functional Analysis and Optimization, 13 (1992), 97-116.

\bibitem {SK}S. Suzuki, D. Kuroiwa, \textit{Fenchel duality for convex set
functions}, Pure and Applied Functional Analysis, 3 (2018), 505--517.

\bibitem {W}D.E. Ward, \emph{Characterizations of strict local minima and
necessary conditions for weak sharp minima}, Journal of Optimization Theory
and Applications, 80 (1994), 551--571.
\end{thebibliography}
\end{document}